\documentclass[12pt]{article}
\usepackage{setspace}
\usepackage{textcomp}
\usepackage[dvips]{color} 
\usepackage{amsmath, amsthm, amssymb}
\usepackage[left=1.4in,
right=1.4in,nohead]{geometry}
\usepackage{hyperref}
\usepackage{latexsym}
\usepackage{mathrsfs}
\usepackage{rawfonts}
\usepackage[dvips]{graphicx}
\usepackage[all]{xy}
\usepackage[sc,osf]{mathpazo}  

\def\bct{\begin{center}}
\def\ect{\end{center}}
\def\beg{\begin}

\def\<{\langle}
\def\>{\rangle}

\def\mbb{\mathbb}
\def\mbbh{\mathbb H}
\def\mbbr{\mathbb R}
\def\mbbz{\mathbb Z}

\def\ni{\noindent}

\def\tn{\textnormal}
\def\wt{\widetilde}

\newtheorem{thm}{Theorem}[section]
\newtheorem{lem}[thm]{Lemma}

\title{Locally conformally flat metrics on surfaces of general type} 
\author{Mustafa Kalafat \and \"Ozg\"ur Kelek\c ci
}

\begin{document}
\maketitle
\begin{abstract} We prove a nonexistence theorem for product type 
manifolds. In particular we show that the 4-manifold  
$\Sigma_g\times\Sigma_h$ does not admit any locally conformally flat  
metric arising from discrete and faithful representations 
for $g\geq 2$ and  $h\geq 1$.

\vspace{.05in}

\ni {\em Keywords:} Hyperbolic geometry; locally conformally flat metric; 4-manifolds.

\vspace{.05in}

\ni {\em Mathematics Subject Classification 2010:} Primary 53C25; Secondary 57M50.
\end{abstract}
\section{Introduction}

A Riemannian $n$-manifold $(M,g)$ is called locally conformally flat (LCF) if any point has a neighborhood which is conformally diffeomorphic to an open subset of $\mbb R^n$. If $n\geq 3$ any conformal map $f: U\longrightarrow \mbb R^n$ from an open subset $U\subset \mbb R^n$ is the restriction of a {\em M\" obius transformation} which is defined to be a composition of inversions on the $n$-sphere. The set of all M\"obius transformations is the conformal transformations of the round sphere. The orientation preserving subset is denoted by $\tn{M\"ob}^+(\mathbb S^n)$ or $\tn{Conf}^+(\mathbb S^n)$.  
In lower dimensions the orientation preserving conformal transformation group can also be realized as the matrix group $\tn{PSL}(2,\mbb F)$ where the field $\mbb F=\mbb R, \mbb C$ or $\mbb H$ for $n=1,2$ or $3$. Conformal diffeomorphisms of the $n$-sphere extends as isometries of the hyperbolic (n+1)-ball that it bounds. Moreover this is a bijective correspondence and we denote the corresponding isometry group by $\tn{Isom}^+\mbbh^{n+1}$. See e.g. \cite{cxslcf} for the details. 
A correct way to understand this isometry group is to realize the hyperbolic space in the Lorentzian $\mbb R^{n+2}$ space as the upper sheet of $-1$ radius hyperboloid. Isometries of the Lorentz metric is
denoted by $\tn{O}(n+1,1)$ which is a subgroup of the general linear group 
$\tn{GL}(n+2,\mbbr)$. The full Lorentz group has $4$ connected components. 
Considering the orientation preserving subset, i.e. of determinant $+1$ 
we get the subgroup $\tn{SO}(n+1,1)$ which has $2$ connected components. 
Among these take the transformations which leaves the upper hyperboloid 
invariant. This is the component which the identity matrix belongs to 
and denoted by $\tn{SO}_e(n+1,1)$. See \cite{hallliegroups} or 
\cite{helgason} for more details on the Lorentz group.  

\vspace{.05in}

Since as a consequence of the Liouville's theorem the transition maps lie 
in a fixed space, one can characterize these spaces in terms of Thurston's 
$(G,X)$ structure terminology in \cite{thurston}.   
One can think of an oriented $n$-manifold with a locally conformally flat 
metric as a $(\tn{Conf}^+\mathbb S^n, \mathbb S^n)$ 
geometry. Moreover, starting from a base point and local conformal chart containing it, considering it as a map $f:U\to \mathbb S^n$ one can embed it into the $n$-sphere. One can embed 
an overlapping chart by pulling back via the transition M\"obius transformations acting on the sphere. Continuing this process one can locally embed i.e. globally immerse the universal covering space $\wt{M}$ of the manifold into the $n$-sphere considering it as a set of classes of loops. Each point of the universal cover represents a loop, and the image of the other endpoint in $\mathbb S^n$ defines the map.  This immersion map 
$D: \wt{M}\longrightarrow \mathbb S^n$ 
is called the {\em developing map} of the LCF manifold. A developing map differs from another by a M\"obius transformation. 
Given any covering transformation $T$ of the universal cover, 
there is a unique element $g\in \tn{Conf}^+\mathbb S^n$ such that 
$D\circ T=g\circ D$. This correspondence defines a homomorphism  
$$\rho : \pi_1(M)\longrightarrow \tn{Conf}^+\mathbb S^n $$
called the {\em holonomy representation} of $M$. 


\vspace{3mm}

Now, we are ready to state the main result of our paper.

\beg{thm} \label{maingroupstthm}
Let $\Gamma=\Gamma_1\times\Gamma_2$ be a product of groups satisfying 
the following three conditions. 
\begin{enumerate}
\item $\Gamma$ is torsion free. 
\item $\Gamma_1$ and $\Gamma_2$ are nontrivial.
\item $\Gamma_1$ is not solvable. 
\end{enumerate}
Then there is no discrete and faithful representation such as $\rho : \Gamma\longrightarrow \tn{SO}_e(n,1).$
\end{thm}

 As an application we can give a partial answer to the open problem 
we described 
in \cite{structuressdlcf}. 

\beg{thm} \label{generaltypesurcaces}
The product 
$\Sigma_g\times\Sigma_h$ of closed surfaces of genus $g\geq 2$ and $h\geq 1$ does 
not admit any locally conformally flat metric 
arising from discrete and faithful representations.
\end{thm}

We can also say something on the possible sign of the locally conformally 
flat metric on these spaces if exists. 

\beg{thm} \label{sign}
The product 
$\Sigma_g\times\Sigma_h$ of closed surfaces of genus $g,h\geq 2$ 
 does not admit any locally conformally flat metric of zero or positive 
 scalar curvature type. Moreover, in the presence of a locally conformally flat metric, the Ricci curvature has to be strictly negative around some point of the manifold. 
\end{thm}

In addition to the references that we mentioned here, the reader is advised to consult to the papers \cite{handle} and \cite{kalafatarguz} for our 
previous work in the field of locally conformally flat geometry.  
In section \S\ref{secrep} we prove the main theorem, and in \S\ref{secscalarcurvature} we analyse the scalar curvature.

\vspace{.05in}

{\bf Acknowledgements.} We thank C. LeBrun, B. Schmidt and J. Morgan for useful suggestions, \c{C}. Karakurt for comments. This work is supported by an RTG grant. Thanks for the Michigan State University for the hospitality.

\section{Representations}\label{secrep}

We start with mentioning a fundamental result, 
{\em Brouwer fixed-point theorem} \cite{hatcherat} which states that any continuous map from the closed $n$-dimensional ball to itself has a fixed point. So any isometry of the hypebolic $n$-space $\mbb H^n$ has at least 
one fixed point either on it or on its ideal boundary $\partial_\infty\mbb H^n$. If you consider the unit ball model, you can think the interior 
as the hyperbolic space and the unit sphere as its ideal boundary at infinity as the following notation suggests.
$$\overline{\mbb H}^n = \mbb H^n \cup \partial_\infty\mbb H^n.$$
A well-known classification of isometries in terms of fixed points can be made as follows, see \cite{benedetti} for a reference. 
If $\gamma\in \tn{Isom}^+\mbb H^n$ is an isometry, then there are 
three mutually exclusive possibilities. 
\beg{enumerate}
\item It has some fixed point(s) in $\mbb H^n$, i.e. in the interior. In this case it is called an {\em elliptic} isometry. 
\item It has exactly one fixed point and that lies at infinity. 
In this case it is called a {\em parabolic} isometry. 
\item It has exactly two fixed points and those lie at infinity. 
In this case it is called a {\em hyperbolic} isometry. 
\end{enumerate}
A nontrivial isometry can not have three fixed points at infinity, since 
then considering the action at boundary it acts on $S^{n-1}$ as 
a M\"obi\"us transformation with three fixed points hence the identity. 
We start with ruling out the first possibility.
\begin{lem}\label{noelliptics} 
If an isometry $\gamma\in\rho(\Gamma)\subset\tn{Isom}^+\mbb H^n$ is 
elliptic and if the subgroup $\<\gamma\>$ it generates is discrete, then 
the order of $\gamma$ is finite. 
\end{lem}
\beg{proof}Under the above assumptions suppose the order of $\gamma$ is infinite. Let $x\in\mbb H^n$ be a fixed point of $\gamma$. 
Then since it is a diffeomorphism $${\gamma_*}_x : T_x\mbbh^n \to T_x\mbbh^n$$ 
is a linear isomorphism preserving the metric, hence an isometry of the 
tangent space at that point. In the following diagram
$$\mbbz\approx \<\gamma\> 
< \tn{Isom}^+_x\mbbh^n 
\stackrel{d_x}{\longrightarrow} \tn{Isom}^+T_x\mbbh^n  
\approx  SO_n$$
the map $d_x$ is injective since in a complete Riemannian manifold, 
an isometry is determined by its value and derivative at a point. 
So that we have an infinite subset of a compact Lie group $SO_n$, 
which should accumulate to say $\varphi\in\tn{Isom}^+\mbbh^n$. 
Suppose $\varphi\notin \<\gamma\>$ otherwise it would obviously be 
contradictory. Now endow the Lie group $\tn{Isom}^+\mbbh^n$ with a left 
invariant Riemannian metric. In this case left multiplication by an 
element is an isometry. Now via discreteness pick an $\epsilon>0$ so that 
$$B_\epsilon(e)\cap \<\gamma\>=\{e\}.$$
Now we claim that this implies 
$$B_\epsilon(\gamma^k)\cap \<\gamma\>=\{\gamma^k\}.$$
Otherwise if there is another element from the subgroup in the ball, 
one can pull this ball back to the one around the identity and find out 
another element there which is pre-avoided.
\end{proof}

Here comes a useful lemma. 
\begin{lem}\label{commutingssame} 
Assume that $\gamma_1,\gamma_2$ are not elliptic and 
$\gamma_1\circ\gamma_2=\gamma_2\circ\gamma_1$,
then they are both parabolic or both hyperbolic with common 
fixed point(s). 
\end{lem}
\beg{proof}Let $\gamma_1$ be a parabolic isometry and 
$x\in\partial_\infty\mbbh^n$ be its unique fixed point. Since
\beg{equation}\label{fixedptsofcommuters} 
\gamma_1(\gamma_2x)=\gamma_2\gamma_1x=\gamma_2x\end{equation}
it fixes $\gamma_2x$ as well so by uniqueness $\gamma_2x=x$. 
So $x$ is also a fixed point of $\gamma_2$ which is on the boundary. 
Then by the trichotomy, $\gamma_2$ can not be elliptic. If $y\neq x$ 
is another fixed point of $\gamma_2$ then  
$$\gamma_2(\gamma_1y)=\gamma_1(\gamma_2y)=\gamma_1y$$
$\gamma_1y$ is also a fixed point of $\gamma_2$ which is distinct from the 
first two, a contradiction. 
So $\gamma_2$ does not have a second fixed point hence it is parabolic as 
well. Replacing the roles of the two isometries, this shows that 
$\gamma_1$ is parabolic if and only if $\gamma_2$ is parabolic.

Secondly suppose $\gamma_1$ is hyperbolic, i.e. two distinct fixed points 
at infinity. Then the commutativity argument (\ref{fixedptsofcommuters}) 
above implies that $\gamma_2$ permutes the fixed points of $\gamma_1$. 
It is easy to see that this permutation is trivial as follows. 
Since there is a unique geodesic joining the fixed points in the 
hyperbolic space (the perpendicular half circle to the boundary plane in 
the upper space model), this geodesic is sent to itself. 
In the ball model this is sending topologically a closed interval to 
itself, so by Brouwer's theorem there is a fixed point along. 
By nonellipticity this fixed point can not lie in the interior hence 
the boundary points are fixed. \end{proof}

Now let $\gamma$ be a parabolic isometry of the hyperbolic $n$-space.  Considering the upper half space model, we can replace $\gamma$ with 
a conjugate so that the fixed point becomes the point $\infty$ at infinity. Then we can write this isometry as follows \cite{benedetti}. 
\beg{equation}\label{parabolicrigidform}
\gamma(x,t)=(Ax+b , t)\end{equation}
for some $A\in SO_{n-1}$ and $b\in \mbbr^n-0$ where $x\in \mbbr^{n-1}$ and 
$t\in\mbbr^+$. This means that $\gamma$ 
acts on the horizontal $(n-1)$-planes by orthogonal transformations. 
If you look at the action on the unit ball this means that $\gamma$ 
leaves the (horospheres) spheres tangent to the boundary at the fixed point. We know that the restriction of the hyperbolic metric to the horospheres or the horizontal planes is a constant multiple of the Euclidean metric, hence this action is isometric on these invariant spaces. We will need the following lemma. 

\begin{lem}\label{nomixedtypeiso} 
Suppose that $\gamma,\bar\gamma$ are two orientation preserving isometries 
of $\mbbh^n$ such that,
\beg{enumerate}
\item $\gamma$ is a parabolic element fixing $x$.
\item $\bar\gamma$ is a hyperbolic element fixing $x$ and some other point $\bar x$.
\end{enumerate}  
Then the subgroup $\<\gamma,\bar\gamma\> < \tn{Isom}^+\mbbh^n$ is not discrete. \end{lem}
\beg{proof}See Lemma D.3.6. of \cite{benedetti}. Working in the upper half space model, we can go through conjugation for $\bar\gamma$ to have fixed points $0,\infty$ and provide $0$ as its common point with $\gamma$. Then we can express $\bar\gamma$ as 
\beg{equation}\label{hyperbolicrigidform} \bar\gamma(x,t)=\lambda(Bx,t)
\end{equation} 
for some  $1\neq\lambda>0$ and $B\in SO_{n-1}$ where $x\in \mbbr^{n-1}$ and $t\in\mbbr^+$. One can think the fixed point at infinity as the point with components $x=0$, $t=\infty$. Now if we conjugate $\gamma$ with $\bar\gamma$ and considering the rigid form (\ref{parabolicrigidform}) of $A$, 
and keeping in mind that the inverse is given by  
$\bar\gamma^{-1}(x,t)=\lambda^{-1}(B^{-1}x,t)$, we have 
{ \renewcommand*{\arraystretch}{1.4}
$$\beg{array}{rcl}
\bar\gamma^{-n} \gamma \bar\gamma^n (x,t) 
&=& \bar\gamma^{-n} \gamma (\lambda^n B^nx,\lambda^n t)\\
&=& \bar\gamma^{-n} (A\lambda^nB^nx+b,\lambda^n t)     \\
&=& (B^{-n}A B^nx+\lambda^{-n}B^{-n}b, t).      
\end{array}$$
}
Now we will assume that $\lambda>1$ otherwise replace $\bar\gamma$ 
with its inverse. Then  
$$\bar\gamma^{-n} \gamma \bar\gamma^n (0_{n-1},1)
=(\lambda^{-n}B^{-n}b, 1)
\longrightarrow (0_{n-1},1) ~~\tn{as}~~n\longrightarrow\infty.$$
The points in this sequence are different from each other and they converge to $(0_{n-1},1)$, consequently the group generated by these 
two transformations does not act properly discontinuously. \end{proof}

Now we have a corollary. 

\beg{thm} If $\Gamma=\Gamma_1\times\Gamma_2$ embeds discretely into the isometry group, then we have exactly one of the two possibilities. Either,
\beg{enumerate}
\item The entire $\rho(\Gamma)$ consists of parabolic elements with common fixed point, or
\item The entire $\rho(\Gamma)$ consists of hyperbolic elements with common fixed points. \end{enumerate} 
\end{thm}
\beg{proof} Because of the Lemma \ref{noelliptics} there are no elliptic elements.  Let $g_1\in \Gamma_1$ be nontrivial and suppose $\rho(g_1,1)$ be parabolic with fixed point $x\in \partial_\infty \mbbh^n$. Then $(1,g_2)$ commutes with that for any $g_2\in \Gamma_2$. So by Lemma \ref{commutingssame}, $\rho(1 \times\Gamma_2)$ is parabolic, 
going back and applying the lemma implies $\rho(\Gamma_1 \times 1)$ 
is also parabolic with the same fixed point. Now taking a general element 
with nontrivial components, 
$$\rho(h_1,h_2)=\rho(h_1,1)\circ \rho(1,h_2)$$
decomposes into parabolic elements that commute hence have the same fixed point. Applying Lemma \ref{commutingssame} again this time to the general element with one of its factors yields that the general element is also parabolic. The hyperbolic case is similar. \end{proof}

Reducing to these two cases, next we will prove that both cases 
are violated. At this point we assume that $\Gamma=\Gamma_1\times\Gamma_2$ 
where $\Gamma_i$ are fundamental groups of hyperbolic surfaces. 
We start with the first possibility. 
\beg{lem} Let $\rho : \Gamma \longrightarrow \tn{Isom}^+\mbbh^n$ be a discrete, faithful representation. Then the image $\rho(\Gamma)$ can not be consist of parabolic isometries with same fixed point. 
\end{lem}
\beg{proof} Consider the upper half space model with entirely parabolic image with common fixed point $\infty$. Then considering the normal form 
\ref{parabolicrigidform} and fixing a horizontal plane (i.e. a horosphere) 
the whole group $\rho(\Gamma)$ acts on a (flat) Euclidean space 
$\mbb E^{n-1}$ again properly discontinuously. Now consider this quotient,
$$M^{n-1}:=\mbb E^{n-1} / \rho(\Gamma)$$
which is a complete, flat Riemannian manifold with $\pi_1M\approx \rho(\Gamma)$. Since the curvature is nonnegative, and the manifold is complete by the Soul theorem \cite{cheegergromollsoul} it has a closed, totally convex, totally geodesic submanifold $S$ called the {\em soul} whose normal bundle is diffeomorphic to the manifold itself. 
By Bieberbach theorem since $S$ is flat it is finitely covered by a torus $\mbb T^k$ so that it has a subgroup of finite index isomorphic to $\mbbz^k$. 
\end{proof}

\ni Next we finally prove the sister theorem for hyperbolic isometries.
\beg{lem} Let $\rho : \Gamma \longrightarrow \tn{Isom}^+\mbbh^n$ be a discrete, faithful representation. Then the image $\rho(\Gamma)$ can not be consist of hyperbolic isometries with same fixed point. 
\end{lem}
\beg{proof} Consider the upper half space model with entirely hyperbolic image with common fixed points $\{0, \infty\}$. For any $\gamma\in\Gamma$ by \ref{hyperbolicrigidform} say $$\rho(\gamma)(x,t)=\lambda_\gamma(A_\gamma x, t)$$ 
for some $\lambda_\gamma>0$ and $A_\gamma\in SO_{n-1}$. Define
$$\varphi: \Gamma \longrightarrow SO_{n-1}\times\mbb R^+
~~~\tn{by}~~~\gamma\mapsto (A_\gamma,\lambda_\gamma).$$
Then for $p_2$ denoting projection onto the second component we have 
$p_2\circ\varphi : \Gamma\to\mbb R^+$. Faithfulness imply that $\varphi$ is injective. Then 
$$\tn{Ker}p_2\varphi=\{\gamma\in\Gamma: \lambda_\gamma=1\}.$$
Since $A_\lambda$ is linear, the isometry $(x,t)\mapsto (A_\lambda x,t)$ 
fixes the whole axis $\{0\}\times \mbb R^+$ hence an elliptic element which is impossible. So this kernel is trivial. $\Gamma$ embeds into the
second component $\mbb R^+$. But it is nonabelian, again a contradiction.
\end{proof}

\section{Scalar curvature}\label{secscalarcurvature}
In this section we will check the possibilities for the sign of the scalar curvature of metrics on the product $4$-manifold. We start with the easiest case. 
\beg{thm} The product $\Sigma_g\times\Sigma_h$ for $g\geq 2$ and $h\geq 1$ does not admit any locally conformally flat metric of zero scalar curvature. 
\end{thm}
\beg{proof} This is a consequence of the Chern-Gauss-Bonnet formula,
\beg{equation}\label{euler} \chi(M)={1\over 8\pi^2}_M\int {s^2\over 24}- 
{|\stackrel{\circ}{\tn{Ric}}|^2 \over 2}+ |W_+|^2+|W_-|^2 \omega_g. \end{equation}
The hypothesis imply that $\chi\leq 0$ and since $\chi=(2-2g)(2-2h)\geq 0$ we have $\chi=0$. So the manifold is Ricci-flat. This means that all the pieces of the curvature tensor vanishes, hence the manifold is flat. 
\end{proof}

\beg{thm} The product $\Sigma_g\times\Sigma_h$ for $g,h\geq 2$ does not admit any metric of positive scalar curvature. 
\end{thm}
\beg{proof}
If an oriented Riemannian manifold has a Riemannian metric of positive scalar curvature, then all of its Seiberg-Witten invariants vanish \cite{witten}.  
So we need to show that Seiberg-Witten invariant of the product 
for $g,h\geq 2$ is non-zero. The Seiberg-Witten invariants are well defined for a compact, oriented 4-manifold with the characteristic number $b_2^+ > 1$ (in our case $b_2^+ = 2gh+1$) \cite{taubesSWGr}. 
Product of K\" ahler metrics is K\"ahler, hence in particular it is symplectic. 
It was shown that  there exists a $ Spin^c $ structure on a closed symplectic 4-manifold for which the associated Seiberg-Witten invariant is 1 \cite{taubesSWSymp}. Thus, the given product manifold can not admit any metric of positive scalar curvature. \end{proof}

\ni Note that the product of hyperbolic metrics on $\Sigma_g\times\Sigma_h$ has negative scalar curvature for $g,h\geq 2$. However the product manifold does not admit any metric of negative sectional curvature for $g,h\geq 1$ because of the 
theorem of Preissmann \cite{byers}. 
On the other hand, the product $\Sigma_g\times\Sigma_1$ for $g\geq 2$ admits metrics of constant scalar curvature of any sign by shrinking or expanding the factors furnished with natural metrics. A result of Schoen and Yau gives an additional information. 

\beg{thm}(\cite{schoenyaukleinian}) 
A compact locally conformally flat Riemannian manifold with $\tn{Ric}_g\geq 0$ 
has universal cover either conformally equivalent to $\mbb S^n$ or isometric to $\mbbr^n, \mbbr\times\mbb S^{n-1}$ where $\mbb S^n$ and $\mbb S^{n-1}$ are spheres of constant curvature. \end{thm}

\ni See \cite{noronha} for a more elementary proof, \cite{zhulcfnonnegativericci} for the generalization to the complete case. In our case for $g,h+1\geq 2$, the universal cover does not include a sphere product, also the Euler characteristic is nonzero so the manifold is not flat. Consequently, this theorem implies that the Ricci curvature has to be strictly negative around some point of the manifolds that we consider.

\bigskip

\bibliography{lcfgeneraltype}{}
\bibliographystyle{alphaurl}


{\small
\beg{flushleft}
\textsc{Orta mh. Z\"ubeyde Han\i m cd. No 5-3 Merkez 74100 Bart\i n, T\" urk\'{i}ye.}\\
\textit{E-mail address:} \texttt{\textbf{kalafat@\,math.msu.edu}}
\end{flushleft}
}

{\small 
\beg{flushleft} \textsc{School of Mathematics, The University of Edinburgh, James Clerk
Maxwell Building, Peter Guthrie Tait Road, Edinburgh EH9 3FD,
Scotland, UK.
}\\
\textit{E-mail address:} \texttt{\textbf{v1okelek@\,ed.ac.uk}}
\end{flushleft}
}

\end{document}